\documentclass[12pt]{article}

\usepackage{latexsym,amsmath,amsthm,amssymb,fullpage}

\usepackage{color}

\usepackage{latexsym,amsmath,amsthm,amssymb,fullpage}

\usepackage[utf8]{inputenc}

\usepackage{amssymb}
\usepackage{epstopdf}
\usepackage{amsmath}
\usepackage{amsthm}
\usepackage{amsfonts}
\usepackage{epstopdf}

\newtheorem{theorem}{Theorem}
\newtheorem{lemma}[theorem]{Lemma}

\newtheorem{corollary}[theorem]{Corollary}

\newtheorem{proposition}[theorem]{Proposition}
\theoremstyle{definition}

\theoremstyle{remark}

\newcommand\beq{\begin {equation}}
\newcommand\eeq{\end {equation}}
\newcommand\beqs{\begin {equation*}}
\newcommand\eeqs{\end {equation*}}

\def\bee{\begin{eqnarray*}}
\def\ene{\end{eqnarray*}}

\newcommand\Z{{\mathbb Z}}
\newcommand\R{{\mathbb R}}
\newcommand\C{{\mathbb C}}
\newcommand\PP{{\mathbb P}}
\newcommand\E{{\mathbb E}}
\newcommand\sph{{\mathbb S}}
\newcommand\T{{\mathbb T}}

\newcommand{\mathbbm}[1]{\text{\usefont{U}{bbm}{m}{n}#1}}

\begin {document}

$$   $$

\vskip 4mm
\begin{center}
\textbf{\LARGE A simple Fourier analytic proof of} \\
\vskip2mm
\textbf{\LARGE the AKT optimal matching theorem} \\ 

\vskip 10mm
{\large Sergey Bobkov and Michel Ledoux}\\
\vskip 3mm
\textit {\large University of Minnesota and University of Toulouse}\\
\end{center}
\vskip 7mm

\vskip 4mm

\begin {abstract}
\noindent
We present a short and elementary proof of the 
Ajtai-Koml\'os-Tusn\'ady (AKT) optimal matching theorem
in dimension 2 via Fourier analysis and a smoothing argument.
The upper bound applies to more general families of samples, including 
dependent variables, of interest in the study of rates of convergence 
for empirical measures. Following the recent pde
approach by L.~Ambrosio, F.~Stra and D.~Trevisan, we also adapt 
a simple proof of the lower bound.
\end {abstract}



\section {Introduction} \label {sec.1}

Given two samples $(X_1, \ldots, X_n)$ and $(Y_1, \ldots, Y_n)$ of 
independent random variables uniformly distributed on the unit square 
$[0,1]^2$, the famous Ajtai-Koml\'os-Tusn\'ady (AKT) optimal matching 
theorem \cite {AKT84} establishes that, with high probability,
$$
\inf_{\sigma}\, \frac 1n \sum_{k=1}^n | X_k - Y_{\sigma(k)}|
	\, \sim \, \sqrt {\frac {\log n}{n} } \, .
$$
Here the infimum is taken over all permutations $\sigma$ of 
$\{1, \ldots, n \}$, $| \cdot |$ is the Euclidean norm in $\R^2$, and 
$ A\sim B$ means that $A \leq C B$ and $B \leq CA$ for some constant 
$C>0$ independent of $n$ $(\geq 2)$. The $\frac 1n$ normalization is 
for convenience with the further statements and formulations, and for 
the purpose of this note, with high probability is simply translated 
by an equivalence on the average
\beq \label {eq.akt}
\E\bigg(\inf_{\sigma}\, \frac 1n \sum_{k=1}^n | X_k - Y_{\sigma(k)}|\bigg)
	\, \sim \, \sqrt {\frac {\log n}{n} } 
\eeq
(concentration arguments allowing for quantitative probabilistic 
estimates cf.~\cite {BL16,AST19}).

The AKT theorem is proved in \cite {AKT84} with combinatorial dyadic
decompositions, where it is also mentioned that the analogous statement 
with the Euclidean norm at the power $p$, $1\leq p < \infty$, holds 
similarly. Further proofs, still based on the same principle, and 
with improved conclusions, have been provided in \cite {S85,S91} or 
\cite {TY93}. M.~Talagrand \cite {T94b,HS92} undertook a deep 
investigation of optimal matching with the tool of the ellipsoid theorem 
from the generic chaining (majorizing measure) theory, with significant 
strengthenings and further, still open, conjectures
(cf.~the monograph \cite {T14}).
In particular, with this approach, he extended in \cite {T92b} 
the upper bound in \eqref {eq.akt} to samples with arbitrary distribution 
on $[0,1]^2$ (which may be then further extended to distributions on 
$\R^2$ under moment conditions \cite {Y92}). Grid matching corresponding 
to $p=\infty$ has been investigated simultaneously \cite {LS89,SY91,T14}...

For the specific uniform distribution, alternate approaches have been 
developed recently, such as gravitational allocation in \cite {HPZ18}. 
A major breakthrough is the investigation
\cite {AST19} by L.~Ambrosio, F.~Stra and D.~Trevisan who used pde methods 
towards exact asymptotics of the optimal matching for $p=2$.

The optimal matching problem may be formulated similarly for samples on 
the cube $[0,1]^d$ for any dimension $d$. The value $d=2$ is actually 
known to be the critical and most delicate one (see the discussion 
in \cite {T14}), since when $d=1$ monotone rearrangement arguments 
show that the order is $\frac {1}{\sqrt n}$ (cf.~\cite {BL16}), while 
when $d \geq 3$, easy tools produce the rate $\frac {1}{n^{1/d}}$, see 
for example \cite {D69} for some early achievements and 
\cite {T94a,DSS13,FG15} for recent more general developments concerning 
$d \geq 3$. We refer in particular to the latter \cite {FG15} and to 
\cite {Y98,T14} for further bibliographical references on the topic of 
optimal matching.

The purpose of this note is to present an elementary Fourier analytic 
proof of the AKT theorem \eqref {eq.akt}, with in particular a very 
simple argument towards the upper bound, valid for any underlying 
distribution on $[0,1]^2$. While the use of Fourier transform is also 
the first step in the Talagrand investigation \cite {T94b,T14} 
(inspired from \cite {CS91}), we replace the delicate genering chaining 
analysis by a standard smoothing procedure. This smoothing procedure is 
also part of the pde analysis developed in \cite {AST19} towards exact 
asymptotics. We borrow from the latter work \cite {AST19} the Lusin 
approximation theorem of Sobolev functions towards a simplified proof 
of the lower bound. The simplicity of the approach developed in this note 
allows for several extensions, and should potentially be useful in the 
study of related issues. Some applications of Fourier analysis and heat 
kernel smoothing in the study of Kantorovich metrics have been proposed 
recently in \cite {S18}.

The note is structured as follows. In Section~\ref {sec.2}, we 
reformulate the optimal matching theorem in suitable Kantorovich metrics 
adapted to Fourier analysis. Next, the main Fourier analytic argument is 
developed, while in Section~\ref {sec.4} the smoothing procedure is 
presented by means of standard Gaussian kernel regularization. The proof 
of the upper bound in the AKT theorem is then immediately deduced in 
Section~\ref {sec.5}, and shown to apply to more general samples, including 
dependent structures. In this formulation, the optimal matching problem 
is part of the study of rates of convergence of empirical measures in 
Kantorovich metrics. Empirical measures with non-random atoms are 
considered in Section~\ref {sec.6}, producing in particular new instances 
of the AKT theorem. The lower bound is established in the next paragraph. 
In the final Section~\ref {sec.8}, we derive more precise quantitative 
upper bounds taking care  of the dependence of the constants as the 
dimension $d$ grows, essentially recovering some claims from \cite {T92a}.

\section {Kantorovich metric} \label {sec.2}

To present the approach, it is convenient to recast the optimal matching 
problem in terms of the Kantorovich metric $\mathrm {W}_1$. 
We mention e.g.~\cite {D02,RR98,V03} as standard references 
on the Kantorovich transport distances.

Given two probability measures $\mu$ and $\nu$ on the Borel sets of $\R^d$
with a first moment, the Kantorovich transport distance $\mathrm {W}_1(\mu,\nu)$ between $\mu$ 
and $\nu$ is defined as
\beq \label {eq.w1}
\mathrm {W}_1(\mu,\nu)
	 \, = \,  \inf_\lambda \int_{\R^d} \int_{\R^d} |x-y| \,d\lambda(x,y)
\eeq
where the infimum is running over all probability measures $\lambda $ on
$\R^d \times \R^d$ with respective marginals $\mu$ and $\nu$, and $|x-y|$ 
represents the Euclidean distance between $x, y \in \R^d$. It is a standard 
consequence of the Birkhoff theorem on the extreme points of the set of 
bi-stochastic matrices that whenever $x_1, \ldots, x_n, y_1, \ldots, y_n$ 
are points in $\R^d$, and $\mu_n = \frac 1n \sum_{k=1}^n \delta_{x_k}$, 
$\nu_n = \frac 1n \sum_{k=1}^n \delta_{y_k}$, then
$$
\mathrm {W}_1(\mu,\nu) \, = \, 
\inf_\sigma \frac 1n \sum_{k=1}^n | x_k - y_{\sigma(k)}| ,
$$
connecting therefore with the optimal matching formulation.
In particular, we will study and state below the AKT results using this
Kantorovich metric $\mathrm {W}_1$.

By the Kantorovich-Rubinstein theorem, the distance $\mathrm {W}_1$
has another description as
\beq \label {eq.kr}
\mathrm {W}_1(\mu,\nu) \, = \ \sup_u \bigg|\int_{\R^d} u \,  d\mu - 
			\int_{\R^d} u \,  d\nu\bigg|
\eeq
where the supremum is taken over all (real-valued) Lipschitz functions $u$ 
on $\R^d$ with Lipschitz semi-norm $\|u\|_{\rm Lip} \leq 1$ with respect 
to the Euclidean distance on $\R^d$.

The aim is to bound the distance $\mathrm {W}_1(\mu,\nu)$ by means of 
Fourier analysis for probability measures supported on a bounded set, 
say $Q^d = (-\pi, \pi]^d$, which requires some additional properties of 
$u$ like periodicity. This is possible, at least when $\mu$ and $\nu$ are 
supported on a smaller part of $Q^d$. In that case, any Lipschitz map $u$ 
on $\R^d$ can indeed be modified outside the sub-cube to become periodic 
and to still be Lipschitz 
(thus not changing the difference of the integrals in \eqref {eq.kr}).

As an alternate approach, one may consider a similar problem on the torus
$\T^d = (\sph^1)^d $ where $\sph^1 = \{z \in \C; \, |z| = 1\}$ denotes the 
unit circle on the complex plane, endowed with the geodesic distance. 
The circle may be identified with the semi-open interval $(-\pi,\pi]$ 
with metric
$$
\rho(x,y) \, = \, \min \big \{|x-y|, 2\pi - |x-y| \big \}, \quad  x,y \in (-\pi,\pi],
$$
via the isometric mapping $U(x) = e^{ix}$. In that case,
$\T^d$ should be identified with $Q^d$ with metric
$$
\rho_d(x,y) \, = \, 
\bigg(\sum_{\ell=1}^d \rho(x_\ell,y_\ell)^2\bigg)^{1/2}, \quad
x = (x_1,\dots,x_d), \ y = (y_1,\dots,y_d) \in Q^d.
$$
For probability measures $\mu$ and $\nu$ on $Q^d$, the Kantorovich 
transport distance with respect to $\rho_d$ is defined similarly to 
\eqref {eq.w1} as
\beq \label {eq.w1tilde}
\widetilde {\mathrm {W}}_1(\mu,\nu) \, = \, 
\inf_\lambda \int_{Q^d} \int_{Q^d} \rho_d(x,y)\,d\lambda (x,y).
\eeq
The general Kantorovich-Rubinstein theorem holds true for the metric 
space $(Q^d,\rho_d)$ as well (cf.~\cite {D02, RR98,V03}) and may be 
restated similarly to \eqref {eq.kr}: For any two Borel probability 
measures $\mu$ and $\nu$ on $Q^d$,
\beq \label {eq.krtilde}
\widetilde {\mathrm {W}}_1(\mu,\nu) 
	\, = \, \sup_u \bigg|\int_{Q^d} u \,d\mu - \int_{Q^d} u \,d\nu\bigg|
\eeq
where the supremum is taken over all (real-valued) maps $u$ on $Q^d$ 
with Lipschitz semi-norm ${\|u\|}_{\textrm {Lip}} \leq 1$ with respect 
to $\rho_d$. The $2\pi$-periodic extention of any such function $u$ 
satisfies $|u(x) - u(y)| \leq {\rm dist}(x-y,2\pi\, \Z^d)$ for all 
$x,y \in \R^d$. In particular, $u$ is continuous on $\R^d$ and has 
Lipschitz semi-norm at most 1 in the sense of the Euclidean distance. 
Conversely, any $2\pi$-periodic function $u$ on $\R^d$ with Euclidean 
Lipschitz semi-norm ${\|u\|}_{\textrm {Lip}} \leq 1$ has 
$\rho_d$-Lipschitz semi-norm at most 1 on $Q^d$. Indeed, using the 
isometric map $U$, it is sufficient to note that the Lipschitz property 
of a function on the torus is a local property, while locally 
the difference between the geodesic and the Euclidean metrics is 
negligible. Thus, the supremum in \eqref {eq.krtilde} may be taken 
over all $2\pi$-periodic $u$ on $\R^d$ with 
${\|u\|}_{\textrm {Lip}} \leq 1$ with respect to the Euclidean distance.

It should also be clear that the supremum in the Kantorovich-Rubinstein representations
may be restricted to $C^\infty$-functions. Once $u$ is $2\pi$-periodic on $\R^d$
and $1$-Lipschitz, i.e. ${\|u\|}_{\textrm {Lip}} \leq 1$, the convolutions
$$
u_\varepsilon(x) \, = \, \frac {1}{(2\pi \varepsilon^2)^{d/2}} \int_{\R^d} 
		u(x - \varepsilon y)\,e^{-|y|^2/2\varepsilon^2}\,dy, \quad x \in \R^d,
$$
of $u$ with Gaussian densities represent $2\pi$-periodic, 
$C^\infty$-smooth, and $1$-Lipschitz functions for any $\varepsilon>0$. 
Since $\max_x |u_\varepsilon(x) - u(x)| \leq d\varepsilon \to 0$ 
as $\varepsilon \to 0$, the function $u$ in \eqref {eq.krtilde}
may be replaced with $u_\varepsilon$'s. A similar remark applies to
the supremum in \eqref {eq.kr} as well.

Since $\rho_d$ is dominated by the usual Euclidean distance, it follows
from \eqref {eq.w1} and \eqref {eq.w1tilde} that
$\widetilde {\mathrm {W}}_1 \leq \mathrm {W}_1$. On the other hand,
$\widetilde {\mathrm {W}}_1(\mu,\nu) = \mathrm {W}_1(\mu,\nu)$ as long 
as both $\mu$ and $\nu$ are supported on a smaller part of the cube 
$Q^d$ such as $[0,\pi]^d$ (suitable for the applications). In this case 
all measures $\lambda $ with marginals $\mu$ and $\nu$ have to be 
supported on $[0,\pi]^d \times [0,\pi]^d$, and since
$\rho_d(x,y) = |x-y|$ in this sub-cube, the right-hand sides 
of \eqref {eq.w1} and \eqref {eq.w1tilde}, and therefore the right-hand 
sides of \eqref {eq.kr} and  \eqref {eq.krtilde}, do coincide.

It is a consequence of this analysis, together with elementary scaling, 
that we may investigate the AKT theorem via the metric 
$\widetilde {\mathrm {W}}_1$ described by \eqref {eq.w1tilde} and 
\eqref {eq.krtilde}. This observation will be used implicitly 
throughout the exposition.

\section{Fourier transform} \label {sec.3}

For a probability measure $\mu$ on the cube $Q^d$, its Fourier-Stieltjes 
transform is defined as the multi-indexed sequence
$$
f_\mu(m) \, = \, \int_{Q^d}  e^{i \left<m,x\right>}\,d\mu(x), \quad  m \in \Z^d,
$$
where $\left<m,x\right> = \sum_{\ell =1}^d m_\ell x_\ell$,
$m = (m_1, \ldots, m_d) \in \Z^d$, $x \in (x_1, \ldots, x_d) \in \R^d$,
which determines $\mu$ in a unique way. Equivalently, $f_\mu$ represents 
the characteristic function of a random vector distributed according to 
$\mu$, which is restricted to the lattice $\Z^d$. Therefore, when bounding 
various distances between two probability measures $\mu$ and $\nu$ on 
$Q^d$, it is sufficient to examine closeness of their Fourier transforms
$f_\mu$ and $f_\nu$.

If a $2\pi$-periodic function $u$ on $\R^d$ is sufficiently smooth,
one may expand it as an absolutely convergent Fourier series
$$
u(x) \, = \, \sum_{m \in \Z^d} a_m\,e^{i\left<m,x\right>}, \quad x \in \R^d,
$$
which can be differentiated term by term. Differentiating this equality 
with respect to the $\ell$-th coordinate, we have
$\partial_\ell u(x) = i\sum_{m \in \Z^d} m_\ell \, a_m \,e^{i\left<m,x\right>}$,
which, according to the Parseval identity, yields
$$
\frac{1}{(2\pi)^d} \int_{Q^d} \big | \partial_\ell u(x) \big|^2 dx \, = \,  
\sum_{m \in \Z^d} m_\ell^2\, |a_m|^2.
$$
Summing over $\ell = 1, \ldots, d$,
$$
\frac{1}{(2\pi)^d} \int_{Q^d} |\nabla u(x)|^2\,dx \, = \, \sum_{m \in \Z^d} |m|^2\, |a_m|^2
$$
where $|m|^2 = \left < m, m\right >$. Moreover, if (additionally) $u$ is $1$-Lipschitz, the modulus
$|\nabla u|$ of its gradient is everywhere less than or equal to $1$, hence
the left-hand side of the preceding is bounded by $1$ so that
\beq \label {eq.lipschitz}
\sum_{m \in \Z^d} |m|^2\, |a_m|^2 \, \leq \, 1.
\eeq

Now, by integration,
$$
\int_{Q^d} u \,d\mu - \int_{Q^d} u \,d\nu  \, = \, 
\sum_{m \not= 0} a_m\,\big [ f_\mu(m) - f_\nu(m) \big].
$$
At this point, the analysis of \cite {T14} makes use of tools from 
the study of stochastic processes. We follow a simpler direct route.
Applying Cauchy's inequality on the basis of \eqref {eq.lipschitz}, 
we arrive at
$$
\bigg|\int_{Q^d} u\,d\mu - \int_{Q^d} u \,d\nu\bigg|^2 \, \leq \,
\sum_{m \not=0} \frac{1}{|m|^2} \, \big |f_\mu(m) - f_\nu(m) \big |^2.
$$
Take then the supremum over all (sufficiently) smooth $2\pi$-periodic 
Lipschitz functions $u$ on the left-hand side with 
${\|u\|}_{\textrm {Lip}} \leq 1$ to reach the following statement.

\begin {lemma} \label {lem.fourier}
Given two probability measures $\mu$ and $\nu$ on
$Q^d$ with Fourier-Stieltjes transforms $f_\mu$ and $f_\nu$, 
$$
\widetilde {\mathrm {W}}_1(\mu,\nu)^2 \, \leq \, 
\sum_{m \not= 0} \frac{1}{|m|^2} \, \big |f_\mu(m) - f_\nu(m) \big |^2.
$$
A similar inequality holds for $\mathrm {W}_1(\mu,\nu)$ 
if $\mu$ and $\nu$ are supported on $[0,\pi]^d$.
\end {lemma}

It could be mentioned that if $\mu$ and $\nu$ have respective smooth
densities $\varphi $ and $\psi$ with respect to $dx$, then
\beq \label {eq.inversesobolev}
\sum_{m \not= 0} \frac{1}{|m|^2} \, \big |f_\mu(m) - f_\nu(m) \big |^2
	\, = \, \frac {1}{(2\pi )^d} \int_{Q^d} 
\big | \nabla \Delta^{-1} (\varphi - \psi) \big|^2 dx
\eeq
where, for a convergent Fourier series 
$g  = \sum_{m \in \Z^d} a_m \, e^{i \left< m, x \right>}$
such that $a_0=0$,
$$
\Delta^{-1} g \, = \,  \sum_{m \not = 0} \frac {1}{|m|^2} \, a_m \, e^{i \left< m, x \right>} .
$$
The quantity on the right-hand side of \eqref {eq.inversesobolev} may be 
identified as an inverse Sobolev-type norm (cf.~\cite {V03}). 
When $\nu = dx$, it has been shown in \cite {P18,S15,L17} that
$$
\widetilde {\mathrm {W}}_2(\mu,\nu)^2 \, = \,  
\inf_\lambda \int_{Q^d} \int_{Q^d} \rho_d(x,y)^2 d\lambda (x,y) \, \leq \,  
\frac {4}{(2\pi )^d} \int_{Q^d} 
\big | \nabla \Delta^{-1} (\varphi - 1) \big|^2 dx .
$$
In this instance ($\nu = dx$), the argument and upper bound developed 
next for the $\widetilde {\mathrm {W}}_1$ distance will therefore apply 
simultaneously to the quadratic Kantorovich distance 
$\widetilde {\mathrm {W}}_2$,
and thus to the AKT theorem \eqref {eq.akt} for $p=2$.

One negative issue about the inequality of Lemma~\ref {lem.fourier} is 
that the sum therein may be divergent. To settle the problem, one may 
use a smoothing operation by suitable convolutions of $\mu$ and $\nu$.

\section{Smoothing} \label {sec.4}

We make use of the simple Gaussian heat kernel smoothing, along the 
line of what is developed in \cite {AST19} (towards more ambitious aims), 
although other convolution kernels might be used to this task. 

On $Q^d$, consider the heat kernel 
$$
p_t (x) \, = \, \frac {1}{(2\pi )^d} 
\sum_{m \in \Z^d} e^{i \left <m , x \right> - |m|^2 t}, 
		\quad t > 0, \,  \, x \in Q^d.
$$
In other words, $p_t$ is the density (with respect to the Lebesgue measure)
of the probability measure $\gamma_t$ supported on $Q_d$ 
whose Fourier-Stieltjes transform is given by
$$
f_{\gamma_t}(m) \, = \, e^{- |m|^2 t},  \quad  m \in \Z^d.
$$
In particular, $\int_{Q^d} p_t(x) dx = 1$. 

If $\mu$ is a probability measure supported on $Q^d$, the heat kernel 
smoothed (probability) measure $\mu_t$, $t >0$, is defined as the convolution 
$\mu * \gamma_t$ via the equality
$$
\int_{Q^d } g \, d\mu_t \, = \,  
\int_{Q^d } \int_{Q^d } g(x-y) p_t(y)\, dy \, d\mu (x)
$$
holding for all $2\pi$-periodic continuous functions $g$ on $\R^d$.
Therefore, if $f_\mu$ is the characteristic function of $\mu$, for every 
$ m \in \Z^d$,
\beq \label {eq.characteristicheat}
f_{\mu_t} (m) \, = \, \int_{Q^d } e^{i \left<m,x\right>} d\mu_t(x) \, = \, 
e^{-|m|^2t} f_\mu (m).
\eeq

The task is now to control the cost in regularization for the Kantorovich metric.
If ${u : \R^d \to \R}$ is $1$-Lipschitz and $2\pi$-periodic, consider
$$
\int_{Q^d} u \, d\mu_t -  \int_{Q^d} u \, d\mu\, = \, 
\int_{Q^d} \int_{Q^d } \big [ u(x-y) - u(x) \big] p_t(y)\, dy\, d\mu(x).
$$
Hence
$$
\bigg |  \int_{Q^d }  u \, d\mu_t -  \int_{Q^d } u \, d\mu \bigg |
	\, \leq \, \int_{Q^d }  |y| \, p_t(y) dy .
$$
and, taking the supremum over all such Lipschitz functions $u$, 
$$
\widetilde {\mathrm {W}}_1 (\mu, \mu_t) \, \leq \, 
\int_{Q^d } |y| \, p_t(y) dy . 
$$

The decay as $t \to 0$ of the expression on the right-hand side actually 
turns out to be of the order of $\sqrt t$. To verify this claim, note that
$\gamma_t$ can be recognized as the product measure whose marginals are 
the image of the Gaussian measure
on the real line with mean zero and variance $2t$ under the map
$$
M(y) \, = \, y - 2\pi k, \quad \pi(2k-1) < y \leq \pi(2k+1), \ k \in \Z.
$$
Indeed, this map pushes forward any probability measure $\eta$ on $\R$
to a probability measure $\widetilde \eta$ on $(-\pi,\pi]$. By the construction,
$M(y) - y$ is a multiple of $2\pi$, so
$ f_{\widetilde \eta}(m) = f_\eta (m)$ for all $ m \in \Z$.
In addition, $|M(y)| \leq |y|$ for all $y \in \R$, so that
$$
\int_{-\infty}^\infty |y|^2\,d\widetilde \eta(y) \, \leq  \, \int_{-\infty}^\infty |y|^2\,d \eta(y).
$$
Choosing for $\eta$ the centered Gaussian measure on the real line 
with variance $2t$, we obtain in this way the one-dimensional marginal 
measure with density $p_t$ on $(-\pi,\pi]$. Moreover, as a consequence 
of the preceding comparison along each coordinate,
$$
\int_{Q^d} |y|^2 p_t(y) dy \, \leq \,  \int_{\R^d} |y|^2 d \eta^{\otimes d}(y) \, = \, 2dt.
$$
As a conclusion of this analysis, for any $\mu$ supported on $Q^d$ 
and any $ t>0$,
\beq \label {eq.regularization}
\widetilde {\mathrm {W}}_1(\mu, \mu_t) \, \leq \,  \sqrt {2dt} \, .
\eeq

We next combine the various steps. By the triangle inequality for 
$\widetilde {\mathrm {W}}_1$ and \eqref {eq.regularization}, for any $t>0$,
$$
\widetilde {\mathrm {W}}_1(\mu,\nu)
	\, \leq \,  \widetilde {\mathrm {W}}_1(\mu_t,\nu_t) +  2 \,\sqrt {2dt} .
$$
It remains to apply the Fourier bound from Lemma~\ref {lem.fourier} to $\mu_t$
and $\nu_t$ which satisfy \eqref {eq.characteristicheat} to reach the following conclusion.

\begin {proposition} \label {prop.main}
Given two probability measures $\mu$ and $\nu$ on
$Q^d$ with Fourier-Stieltjes transforms $f_\mu$ and $f_\nu$, for any $ t >0$,
$$
\widetilde {\mathrm {W}}_1(\mu,\nu)
	\, \leq \, \bigg (\sum_{m \not= 0} \frac {1}{|m|^2} \, e^{-2 |m|^2t} \,
	 \big | f_\mu (m) - f_\nu (m) \big|^2 \bigg)^{1/2}
    		+  2 \,\sqrt {2dt} .
$$
A similar inequality holds for $\mathrm {W}_1(\mu,\nu)$ if $\mu$ and $\nu$ 
are supported on $[0,\pi]^d$.
\end {proposition}

\section {Application to the AKT theorem} \label {sec.5}

This section describes the application of the preceding Fourier analytic 
approach to the upper bound in the AKT theorem. It actually applies to 
a somewhat extended probabilistic setting, a form of which having 
already been emphasized in \cite {T92b}.

Namely, consider random variables $X_1, \ldots, X_n,Y_1, \ldots, Y_n$ 
with values in $[0,1]^d$ such that the couples $(X_1,Y_1), \ldots, (X_n,Y_n)$ 
are pairwise independent and, for every $k=1, \ldots, n$, $X_k$ and $Y_k$ 
have the same distribution. Apply Proposition~\ref {prop.main} to 
the empirical measures
$\mu_n = \frac{1}{n}\,\sum_{k=1}^n \delta_{X_k}$ and
$\nu_n = \frac{1}{n}\,\sum_{k=1}^n \delta_{Y_k}$ 
(supported on $[0,1]^d \subset [0,\pi]^d$) to get that, 
after averaging and use of Jensen's inequality,
$$
\E \big ({\mathrm {W}}_1(\mu_n,\nu_n) \big)
	\, \leq \, \bigg (\sum_{m \not= 0} \frac {1}{|m|^2} \, e^{-2|m|^2t} \,
	\E \big ( \big | f_{\mu_n} (m) - f_{\nu_n} (m) \big|^2 \big)\bigg)^{1/2}
    		+  2\, \sqrt {2dt} 
$$
for any $t >0$. Now, by the independence and equidistribution assumptions 
on the variables $X_1, \ldots, X_n,Y_1, \ldots, Y_n$,
$$
\E \big ( \big | f_{\mu_n} (m) - f_{\nu_n} (m) \big|^2 \big) \, \leq \,  \frac 4n
$$
for every $m \in \Z^d$, so that
\beq \label {eq.optimt}
\E \big ({\mathrm {W}}_1(\mu_n,\nu_n) \big)
	\, \leq \, \frac {2}{\sqrt n} \bigg (\sum_{m \not=0} 
\frac {1}{|m|^2} \, e^{-2|m|^2t} \bigg)^{1/2}
    		+  2\, \sqrt {2dt} .
\eeq
From a (crude) comparison between series and integral,
it should be clear without computations that, up to $d$-dependent factors,
\beq \label {eq.sd}
S_d(t) \, = \, \sum_{m \not= 0} \frac {1}{|m|^2} \, e^{-2|m|^2t}
  \, \sim \, \int_{|x| \geq 1} \frac{1}{|x|^2}\, e^{-2t |x|^2}\,dx 
  \, \sim \, \int_1^\infty r^{d-3}\, e^{-2tr^2} dr .
\eeq
For the small values of $t>0$, the latter integral is of order $1$ if $d=1$, 
$\log (\frac 1t)$ if $d=2$ and $t^{-(d/2) +1}$ if $d \geq 3$.
After optimization in $t >0$ in \eqref {eq.optimt}, we thus conclude to 
the following statement which covers the upper bound in the AKT theorem 
when $d=2$,
providing at the same time the optimal rates for $d=1$ and $d \geq 3$.

\begin {theorem} \label {thm.main}
Let $X_1, \ldots, X_n,Y_1, \ldots, Y_n$ be random variables with 
values in $[0,1]^d$ such that the couples $(X_1,Y_1), \ldots, (X_n,Y_n)$ 
are pairwise independent and, for every $k=1, \ldots, n$, 
$X_k$ and $Y_k$ have the same distribution. For the empirical measures
$\mu_n = \frac 1n \sum_{k=1}^n \delta_{X_k}$ and 
$\nu_n = \frac 1n \sum_{k=1}^n \delta_{Y_k}$ associated to the samples 
$(X_1, \ldots, X_n)$ and $(Y_1, \ldots, Y_n)$, it holds true that
\beqs
 \E \big ( \mathrm {W}_1 (\mu_n, \nu_n) \big)  \, = \,
 \begin {cases}
  O \big (\frac {1}{\sqrt n} \big)  & \text{if \, $d=1$,} \\
    O \Big (\sqrt {\frac {\log n}{n}} \, \Big)  & \text{if \, $d=2$,} \\
   O \big (\frac {1}{n^{1/d}} \big)  & \text {if \, $ d \geq 3$.} \\
\end {cases}
\eeqs
\end {theorem}

In the last Section~\ref {sec.8}, we develop a more careful analysis 
of the function $S_d(t)$ of \eqref {eq.sd} to reach more explicit 
quantitative bounds, in particular with respect to dependence as 
the dimension $d $ increases. Namely, Proposition~\ref {prop.quantitative} 
below with $\delta = \frac {2}{\sqrt n}$ yields the following quantitative 
statement of Theorem~\ref {thm.main},
\beq \label {eq.quantitative}
 \E \big ( \mathrm {W}_1 (\mu_n, \nu_n) \big)  \, \leq \,
 \begin {cases}
 \frac {2}{\sqrt n}   & \text{if \, $d=1$,} \\
     10 \, \sqrt {\frac {1 + \log n}{n}}  & \text{if \, $d=2$,} \\
   \frac {16 \sqrt {d}}{n^{1/d}}  & \text {if \, $ d \geq 3$.} \\
\end {cases}
\eeq
The numerical constants are not sharp, but the order of growth as 
$d \to \infty$ matches the first order asymptotics of \cite {T92a}.

If the random variables $X_1, \ldots, X_n, Y_1, \ldots, Y_n$ are independent
and have the same law $\mu$, then by Jensen's inequality
$$
\E \big ( \mathrm {W}_1 (\mu_n, \nu_n) \big) \, \geq \, 
	\E \big ( \mathrm {W}_1 (\mu_n, \mu) \big)
$$
since $\E(\nu_n) = \mu$. The upper bounds of Theorem~\ref {thm.main} 
and \eqref {eq.quantitative} thus apply to $\E  (\mathrm {W}_1 (\mu_n, \mu))$.
As such, the conclusions enter the framework of rates of convergence
for empirical measures.

As an example illustrating Theorem 3, one may consider two sequences
$$
	X_k(\omega) \,= \, U(k \omega_1 + \omega_2), \quad 
	Y_k(\omega) \, =\,  V(k \omega_1 + \omega_2), \quad 
	\omega = (\omega_1,\omega_2) \in \Omega, \quad k \geq 1,
$$
defined for given Borel measurable functions $U,V:[0,1] \to [0,1]^d$
on the square $\Omega = [0,1] \times [0,1]$, which we equip with 
the normalized Lebesgue measure $\mathbb P$. As easy to check, 
${(X_k)}_{k \geq 1}$ forms a strictly stationary sequence of pairwise 
independent random variables (which however are not independent), and 
the same is true for ${(Y_k)}_{k \geq 1}$. If $U$ and $V$ have equal 
distributions under the Lebesgue measure on $[0,1]$, then 
Theorem~\ref {thm.main} is applicable, so that one can make 
the conclusion about the closeness of the associated empirical measures.

One may even further generalize Theorem~\ref {thm.main} to the setting 
of weakly dependent random variables. Recall that, given a probability 
space $(\Omega,\mathcal {A},\mathbb P)$ and two $\sigma$-algebras 
$\mathcal {A}_1, \mathcal {A}_2 \subset \mathcal {A}$, the Rosenblatt
coefficient, which quantifies the strength of dependence between 
${\cal A}_1$ and ${\cal A}_2$, is defined to be
$$
\alpha(\mathcal {A}_1, \mathcal {A}_2) \, = \,
	\sup\big\{ \big |\mathbb P(A_1 \cap A_2) - \mathbb P(A_1) \, \mathbb P(A_2) \big | ; \, 
	A_1 \in \mathcal {A}_1,  A_2 \in \mathcal {A}_2\big\}.
$$
It is one of eight well known measures of dependence (and the weakest one) 
which is used in the theory of strong mixing conditions 
(cf.~\cite{B05,B07}). Clearly,
$$
\alpha(\mathcal {A}_1, \mathcal {A}_2) \, = \, \sup \big |\mathrm {Cov}(\varphi,\psi) \big |
$$
where the supremum is running over all $\mathcal {A}_1$- and 
respectively $\mathcal {A}_2$-measurable functions $\varphi$ and $\psi$ 
on $\Omega$ with values in $[0,1]$. If $\varphi$ and $\psi$ are 
complex-valued with $|\varphi| \leq 1$ and $|\psi| \leq 1$, then,
by the bilinearity of the covariance functional,
$\mathrm {Cov}(\varphi,\psi) = 
\E \big ((\varphi - \E (\varphi)) (\overline \psi - \E (\overline \psi))\big )$ 
is bounded in absolute value by $16\,\alpha(\mathcal {A}_1, \mathcal {A}_2)$.

In practice, one is given a sequence of $\sigma$-algebras $\mathcal {A}_k$ 
generated by random elements $Z_k$, $k \geq 1$, defined on the same 
probability space $\Omega$, with which one associates the characteristics
$$
\alpha(\ell) \, = \,  \sup_{|j - k| \geq \ell} \alpha(\mathcal {A}_j, \mathcal {A}_k), \quad  \ell \geq 1.
$$
Repeating the arguments in the proof of Theorem~\ref {thm.main}, we have:

\begin {corollary}
Let $Z_k = (X_k,Y_k)$, $k \geq 1$, be random variables with values in 
$[0,1]^d \times [0,1]^d$ such that $X_k$ and $Y_k$ have the same distribution 
for every $k \geq 1$. If the associated mixing sequence 
$\alpha(\ell)$, $\ell \geq 1$, is summable, then the asymptotic bounds 
of Theorem~\ref {thm.main} remain to hold for the empirical measures
$\mu_n = \frac 1n \sum_{k=1}^n \delta_{X_k}$ and 
$\nu_n = \frac 1n \sum_{k=1}^n \delta_{Y_k}$.
\end {corollary}

Indeed, for every $m \in \Z^d$,
\beqs \begin {split}
\E \big ( \big | f_{\mu_n} (m) - f_{\nu_n} (m) \big|^2 \big) 
	& \, = \,  
\frac{1}{n^2}\, \sum_{j,k=1}^n 
\mathrm {Cov} \big(e^{i\left<m,X_j\right>} - e^{i\left<m,Y_j\right>},
e^{i\left<m,X_k\right>} - e^{i\left<m,Y_k\right>}\big)  \\
	& \, \leq \, 
\frac 4n + \frac{64}{n^2} \sum_{1 \leq j \neq k \leq n} 
\alpha \big (|j - k| \big ) \\
	& \, = \, 
\frac 4n + \frac{128}{n^2} \, 
\sum_{\ell =1}^{n-1} \, (n- \ell)\, \alpha(\ell)  \\
	& \, \leq \, 
\frac 4n + \frac{128}{n} \, \sum_{\ell =1}^{n-1} \, \alpha(\ell).
\end {split} \eeqs
It therefore remains to apply Proposition 2 as for Theorem~\ref {thm.main}.

\section {Empirical measures with non-random atoms} \label {sec.6}

As another application of the preceding approach, fix a collection 
of points in the unit cube $[0,1]^d$, say $x_1,\dots,x_N$, $N \geq 2$.
One may use various selections of indices to construct (deterministic) 
empirical measures with atoms at $x_j$ (repetition of the points 
in the sequence is allowed). Namely, for $1 \leq n \leq N$,
let $\mathcal {G}_n$ denote  the collection of all subsets $\tau$ 
of $\{1,\dots,N\}$ of cardinality $|\tau| = n$ equipped with 
the uniform probability measure $\pi_n$. With every 
$\tau \in \mathcal {G}_n$, we associate an ``empirical" measure
$$
\mu_\tau \, = \,  \frac{1}{n} \sum_{j \in \tau} \delta_{x_j},
$$
which may be treated as a random measure on the probability space
$(\mathcal {G}_n,\pi_n)$. The goal is to show that most of 
$\mu_{\tau}$'s are concentrated around the average measure
\beq \label {eq.mutau}
\mu \, = \, \E_{\pi_n} (\mu_\tau) \, = \, 
\int_{\mathcal {G}_n} \mu_\tau\, d\pi_n(\tau) \, = \, 
\frac{1}{N} \sum_{j= 1}^{N} \delta_{x_j}
\eeq
as long as $n$ is large (in the sense of the distance $\mathrm {W}_1$). 
For simplicity, we skip the parameter $N$ since interest in the final 
estimates is concerned with the dependence with respect to the growing 
$n$, while $N$ may be arbitrarily large.
To this aim, consider the functional
$$
L_u(\tau) \, = \,  \int_{\mathcal {G}_n} u(x)\,d\mu_\tau(x) 
	\, = \,  \frac{1}{n} \sum_{j \in \tau} u(x_j),  \quad \tau \in \mathcal {G}_n,
$$
associated to a given complex-valued function $u$ on the cube $[0,1]^d$. 
As is easy to check,
$$
\mathrm {Var}_{\pi_n}(L) \, = \, 
\E_{\pi_n}\big ( \big | L - \E_{\pi_n} (L) \big |^2\big) \, = \, 
\frac{N-n}{2n N^2(N-1)} \, \sum_{i,j = 1}^N \big |u(x_i) - u(x_j) \big |^2.
$$
If $|u| \leq 1$, it follows that 
$\mathrm {Var}_{\pi_n}(L_u)  \leq  \frac{2}{n}$. Since the 
Fourier-Stieltjes transform $f_{\mu_\tau}( \pi m)$ corresponds to 
$L_u(\tau)$ with $u(x) = e^{i \pi \left<m,x\right>}$, the analysis of 
the preceding section may be developed in the same way. Together with 
the more quantitative estimates from Proposition~\ref {prop.quantitative} 
below, the following corollary holds true.

\begin {corollary} \label {cor.tau}
Given a collection of points $x_1,\dots,x_N$ in $[0,1]^d$,
for any integer $1 \leq n \leq N$, the empirical measures $\mu_\tau$ satisfy
\beqs
 \E_{\pi_n} \big ( \mathrm {W}_1 (\mu_\tau, \mu) \big)  \, \leq \,
 \begin {cases}
  \sqrt {\frac 2n}   & \text{if \, $d=1$,} \\
    8 \, \sqrt {\frac {1 + \log(2n)}{n}}   & \text{if \, $d=2$,} \\
  \frac {13\sqrt {d}}{n^{1/d}}  & \text {if \, $ d \geq 3$.} \\
\end {cases}
\eeqs
\end {corollary}

Note in particular that if $N = 2n$,
$$
\E_{\pi_n} \big ( \mathrm {W}_1 (\mu_\tau, \mu) \big) \, = \,  
\frac{1}{2n}\, \E \bigg ( \inf \ \sum_{k=1}^n |x_{i_k} - x_{j_k}|  \bigg)
$$
where the averaging on the right is performed over all choices of indices
$1 \leq i_1 < \dots < i_n \leq 2n$, while the infimum is taken over
all permutations $j_1,\dots,j_n$ of the remaining integers in the set
$\{1,\dots,2n\} \setminus \{i_1,\dots,i_n\}$.

In fact, the preceding corollary easily implies Theorem~\ref {thm.main} 
specialized to the iid case. This is achieved by averaging 
\eqref {eq.mutau} over $x_1,\dots,x_{2n}$ according to the product measure
$\mu^{\otimes 2n}$ for a fixed probability distribution $\mu$ on $[0,1]^d$.
Actually, the argument extends to more general classes. Namely, if the
joint distribution of the random vectors $X_1,\dots,X_{2n}$ with values 
in $[0,1]^d$ is invariant under permutations of the indices, then for 
the empirical measures $\mu_n = \frac{1}{n}\sum_{k=1}^n \delta_{X_k}$
and $\mu$ the distribution of $X_1$, $\E (\mathrm {W}_1 (\mu_n, \mu))$ 
is controlled as in Corollary~\ref {cor.tau}.

\section {Lower bound} \label {sec.7}

While the AKT upper bound may be extended to families of samples with 
arbitrary (compactly supported) distributions, it is well-known 
(cf.~e.g.~\cite {BB13,T14}) that the lower bound requires distributions 
with enough regularity, for example absolutely continuous with respect 
to Lebesgue measure. A pde proof of the lower bound in the AKT theorem 
has been provided recently in the paper \cite {AST19}, relying on a somewhat heavy
analysis involving in particular Riesz transform bounds. We extract here 
the necessary argument in our framework via a simple fourth moment
computation, thereby producing a rather mild proof.

Let $X_1, \ldots, X_n, Y_1, \ldots, Y_n$ be independent with uniform
distribution $d\mu = \frac {dx}{(2\pi)^d} $ on $Q^d$.
For any $t>0$, contractivity of the Kantorovich metric shows that
$$
\widetilde {\mathrm {W}}_1(\mu_n, \nu_n) \, \geq \, 
	\widetilde {\mathrm {W}}_1(\mu_{n,t}, \nu_{n,t}).
$$
This is actually immediate from the definition of 
$\widetilde {\mathrm {W}}_1$ 
and the heat kernel regularization since
$$
\int_{Q^d} u \, d\mu_{n,t} - \int_{Q^d} u \, d\mu
	\, = \, \int_{Q^d} \bigg [ \int_{Q^d} u(x-y)\, d\mu_n(x) 
		-  \int_{Q^d} u(x-y)\, d\mu(x) \bigg] p_t(y)\, dy 
$$
and $\int_{Q^d} p_t(y)\, dy = 1$.

Using an absolutely convergent random Fourier series, let
$$
h(x) \, = \,  \sum_{m \not= 0} \frac {1}{|m|^2} \, e^{-|m|^2 t} 
\bigg ( \frac 1n \sum_{k=1}^n 
\big [ e^{i \left<m,X_k\right>} - e^{i \left<m,Y_k\right>}\big]\bigg)\, 
e^{- i \left<m,x\right>} , \quad x \in \R^d.
$$
This equality defines a $2\pi$-periodic, real-valued, $C^\infty$-smooth function,
whose Laplacian 
$$
\Delta h(x) \, = \, - 
\sum_{m \not= 0} e^{-|m|^2 t} \bigg ( \frac 1n \sum_{k=1}^n 
	\big [ e^{i \left<m,X_k\right>} - e^{i \left<m,Y_k\right>}\big]\bigg)\,  e^{- i \left<m,x\right>}
$$
represents the multiple Fourier series for the density of
$\nu_{n,t} - \mu_{n,t}$ (with respect to the Lebesgue measure on $Q^d$). 
Hence, the integration by parts formula for a smooth $2\pi$-periodic function 
$v : \R^d \to \R$ yields
\beq \label {eq.ip}
\int_{Q^d} v \, d\mu_{n,t} - \int_{Q^d} v \, d\nu_{n,t} \, = \, - 
\int_{Q^d} v \, \Delta h\,d\mu 
	\, = \,   \int_{Q^d} \left<\nabla h,\nabla v \right> d\mu.
\eeq

For $\alpha >0$, denote by $u : Q^d \to \R$ the $\alpha$-Lipschitz 
(i.e. ${\| u \|}_{\textrm {Lip}} \leq \alpha$)
Lusin extension of $h$ on the torus $(Q^d, \rho_d)$ such that
\beq \label {eq.lusin}
\mu \big ( \{ h \not= u \} \big) \, \leq \, 
\frac {K}{\alpha^2} \int_{Q^d} |\nabla h |^2 d\mu
\eeq
where $K > 0$ only depends on $d$ 
(\cite {AF84}, cf.~Lemma 5.1 in \cite {AST19}). 
By the Kantorovich-Rubinstein theorem,
$$
\widetilde {\mathrm {W}}_1(\mu_n, \nu_n) \, \geq \,  
\widetilde {\mathrm {W}}_1 (\mu_{n,t} , \nu_{n,t}) \, \geq \, 
	\frac {1}{\alpha} \,  
\bigg|\int_{Q^d} u \,d\mu_{n,t} - \int_{Q^d} u \,d\nu_{n,t}\bigg| .
$$
On the other hand, by \eqref {eq.ip},
\beqs \begin {split}
	\int_{Q^d} u \, d\mu_{n,t} - \int_{Q^d} u \, d\nu_{n,t} 
	& \, = \,   \int_{Q^d} \left<  {\nabla h}, \nabla u \right> d\mu \\
	& \, = \,  \int_{Q^d} |\nabla h|^2 \, d\mu -
		\int_E \left<\nabla h,  \nabla h - \nabla u \right> d\mu \\
\end {split} \eeqs
where $E = \{ h \not= u \}$. Hence
$$
 \alpha \, \E \big ( \, \widetilde {\mathrm {W}}_1(\mu_n, \nu_n) \big)  \, \geq \, 
\E \bigg ( \int_{Q^d} |\nabla h|^2 \, d\mu \bigg) - \E \bigg 
( \bigg |\int_E \left<\nabla h, \nabla h - \nabla u \right> d\mu \bigg| \bigg). 
$$

Since the $X_k,Y_k$'s are independent and uniformly distributed on $Q^d$,
\begin {equation*} \begin {split}
\E \bigg ( \int_{Q^d} |\nabla h|^2 \, d\mu \bigg) 
	 & \, = \, \sum_{m\not=0} \frac {1}{|m|^2} \, e^{-2|m|^2t} \,
	 		\frac {1}{n^2} \,\E \bigg ( \bigg | \sum_{k=1}^n
\big[e^{i \left<m,X_k\right>} - e^{i \left<m,Y_k\right>}\big]\bigg|^2\bigg) \\
	 & \, = \, \frac 1n \sum_{m\not=0} \frac {2}{|m|^2} \, e^{-2|m|^2t} .
\end {split} \end {equation*}
Denote by $c(n,t)$ this quantity, where $t = t(n) \in (0,1)$ will be specified.
On the other hand, since $u$ is $\alpha$-Lipschitz, and by repeated use of H\"older's inequality,
\begin {equation*} \begin {split}
\bigg | \int_E \left<\nabla h,\nabla h - \nabla u \right> d\mu \bigg |
	& \, \leq \, \int_E |\nabla h|^2 \, d\mu + \alpha \int_E |\nabla h| \, d\mu \\
	& \, \leq \, \mu (E)^{1/2} \bigg (\int_{Q^d} |\nabla h|^4 \, d\mu \bigg)^{1/2}
	+ \alpha \mu (E)^{3/4} \bigg (\int_{Q^d} |\nabla h|^4 \, d\mu \bigg)^{1/4}
\end {split} \end {equation*}
and
\begin {equation*} \begin {split}
\E \bigg (& \bigg | \int_E 
    \left<\nabla h, \nabla h - \nabla u \right> d\mu \bigg | \bigg) \\
	&  \, \leq \, \big [\E \big(\mu (E)\big) \big] ^{1/2}  
	\, \bigg [\E \bigg(\int_{Q^d} |\nabla h|^4 \, d\mu \bigg) \bigg]^{1/2}  
		+ \alpha \, \big [ \E \big(\mu (E)\big) \big]^{3/4}  \, 
	\bigg [ \bigg (\E \int_{Q^d} |\nabla h|^4 \, d\mu \bigg) \bigg]^{1/4} . \\
\end {split} \end {equation*}
Moreover, by the Lusin approximation \eqref {eq.lusin},
$$
\E  \big (\mu (E) \big) \, \leq \, \frac {K}{\alpha^2} \, c(n,t) .
$$ 
If we let $d(n,t) = \E \big (\int_{Q^d} |\nabla h|^4 \, d\mu \big)$, 
we have therefore obtained that
\beq \begin {split} \label {eq.lusincd}
 \alpha \, \E \big ( \, \widetilde {\mathrm {W}}_1(\mu_n, \mu) \big) 
		& \, \geq \,  c(n,t) - \frac {1}{\alpha} \, \big (Kc(n,t) \big)^{1/2}\, d(n,t)^{1/2} \\
		& \quad \, \,	-  
		\frac {1}{\sqrt {\alpha}} \, \big (Kc(n,t) \big)^{3/4}\, d(n,t)^{1/4} . \\
\end {split} \eeq

Before optimization of the choice of $\alpha >0$, we need to evaluate 
$d(n,t)$. By the triangle inequality,
$$
d(n,t) \, \leq \, 8 \ \E \bigg (\int_{Q^d} |\nabla \widetilde {h}|^4 d\mu\bigg)
$$
where $\widetilde {h}(x) = \sum_m b_m e^{-i \left<m,x\right>}$ with
$$
b_m \, = \,  \frac {1}{|m|^2} \, e^{-|m|^2t}  
\bigg ( \frac 1n \sum_{k=1}^n e^{i \left<m,X_k\right>} \bigg)
$$
for $m \not=0$ and $b_0 = 0$. It holds that
$$
\big |\nabla \widetilde {h}(x) \big |^2  \, = \, - \sum_{m_1,m_2 \in \Z^d} 
 	\left<m_1,m_2\right> b_{m_1} b_{m_2} \,e^{- i \left<m_1 + m_2,x\right>}
$$
and
$$
\int_{Q^d} |\nabla \widetilde {h}|^4\, d\mu  \, = \, 
   \sum \left<m_1,m_2\right> b_{m_1} b_{m_2} \left<m_3,m_4\right> b_{m_3} b_{m_4}
$$
where the sum is taken over $m_1,m_2,m_3,m_4 \in \Z^d $ such that 
$m_1 + m_2 + m_3 + m_4 = 0$.
Now
\beqs \begin {split}
\E (  b_{m_1} b_{m_2} \, b_{m_3} b_{m_4}) 
   &\, = \, \frac {1}{|m_1|^2|m_2|^2|m_3|^2|m_4|^2} \,
  e^{-(|m_1|^2 + |m_2|^2 + |m_3|^2 + |m_4|^2)t} \\
 &\quad \, \, \cdot \frac {1}{n^4} \sum_{k_1,k_2,k_3,k_4=1}^n 
\E \big(e^{i \langle m_1,X_{k_1}\rangle} \, e^{i \langle m_2,X_{k_2}\rangle} \,
e^{i \langle m_3,X_{k_3}\rangle} \,e^{i \langle m_4,X_{k_4}\rangle}\big) . \\
\end {split} \eeqs
Since the relevant indices satisfy $m_\ell \not= 0$, $\ell = 1,2,3,4$, and 
$m_1 + m_2 + m_3 + m_4 = 0$, the last expectation
is non zero, equal to $1$, only if $k_1 = k_2 = k_3 = k_4$ or if
$$
\begin {cases}
k_1=k_2, \, \, k_3=k_4  & \text{and} \quad m_1 + m_2 =0, \, \, m_3 + m_4=0, \\
k_1=k_3,\, \,  k_2=k_4  & \text{and} \quad m_1 + m_3 =0, \, \, m_2 + m_4 =0, \\
k_1=k_4, \, \, k_2=k_3  & \text{and} \quad m_1+ m_4 =0, \, \, m_2 +m_3 =0. \\
\end {cases}
$$
These respective contributions yield the upper bound
\begin {equation*} \begin {split}
\E \bigg ( \int_{Q^d} |\nabla \widetilde {h}|^4\, d\mu \bigg) 
	&\, \leq \, \frac {1}{n^3} \, \sum
	\frac {1}{|m_1| |m_2| |m_3| |m_4|} 
			\,  e^{-(|m_1|^2 + |m_2|^2 + |m_3|^2 + |m_4|^2)\,t} \\
	& \quad \, \, + \frac {3}{n^2} \sum_{m_1 ,m_3 \not= 0}
	\frac {1}{|m_1|^2 |m_3|^2} \, e^{-2(|m_1|^2 + |m_3|^2)\,t}  \\
	& \, = \, e(n,t) + \frac 34 \, c(n,t)^2 
\end {split} \end {equation*}
where the first sum on the right-hand side is over all 
$m_1,m_2,m_3,m_4 \in \Z^d \setminus \{0\}$. 

It is easily seen that, for $0 < t \leq \frac 12$,
$$
e(n,t) \, = \, \frac {1}{n^3} 
\bigg ( \sum_{m \not= 0} \frac {1}{|m|} \,  e^{-|m|^2 t}  \bigg)^4
 \, \sim \, \frac {1}{n^3} \, \frac {1}{t^{(d-1)/2}} 
\int_{\sqrt t}^\infty r^{d-2}e^{-r^2} dr
$$
while
$$
 c(n,t) \, \sim \, \frac {1}{n} \, \frac {1}{t^{(d-2)/2}} 
\int_{\sqrt {2t}}^\infty r^{d-3}e^{-r^2} dr. 
$$
In the following, take $d=2$.  Hence $e(n,t)$ is of the order of 
$\frac {1}{n^3\sqrt t} $
and $c(n,t)$ of the order of $\frac 1n \log (\frac 1t)$.
Choosing $t = t(n) = \frac {1}{2n}$, we see that
$e(n,t)$ is negligeable with respect to the square of
$$
	c_n \, = \, c\big( n, t(n) \big )  \sim \,  \frac {\log n}{n} \, .
$$
So for this choice of $t=t(n)$, for some constant $K' >0$, 
$$
 d\big( n, t(n) \big) \, = \, 
\E \bigg ( \int_{Q^d} |\nabla h|^4\, d\mu \bigg) \, \leq \,  K' c_n^2.
$$

Summarizing these estimates in \eqref {eq.lusincd} for 
$t = t(n) = \frac {1}{2n}$, we then get
$$
	\alpha \, \E \big ( \, \widetilde {\mathrm {W}}_1(\mu_n, \nu_n) \big)  
	\, \geq \, c_n - \frac {1}{\alpha} \,  ( K c_n )^{1/2}
( K' c_n^2)^{1/2} - \frac {1}{\sqrt \alpha} \, ( K c_n)^{3/4}(K'c_n^2)^{1/4} .
$$
For the choice of $\alpha = \beta \sqrt {c_n}$ with $\beta >0$ 
large enough, it follows that
$$
\E \big ( \, \widetilde {\mathrm {W}}_1(\mu_n, \nu_n) \big)  
		\, \geq \, c \sqrt {c_n}  \, \sim \, \sqrt { \frac {\log n}{n} }
$$
which is the expected lower bound in the AKT theorem \eqref {eq.akt}.

It should be mentioned that in dimension one, the lower bound of the order
of $\frac {1}{\sqrt n}$ is easily achieved via the monotone representation 
$$
\mathrm {W}_1 (\mu_n , \nu_n) \, = \,  
\int_0^1 \frac 1n \,  \bigg | \sum_{k=1}^n 
\big(\mathbbm{1}_{\{X_k \leq x\}} - \mathbbm {1}_{\{Y_k \leq x\}}\big)\bigg|\,dx
$$
for independent uniform random variables 
$X_1, \ldots, X_n, Y_1, \ldots, Y_n$ on $[0,1]$
(cf.~\cite {BL16}). Hence
$$
\E \big ( \mathrm {W}_1 (\mu_n , \nu_n) \big ) \, \geq \, 
\frac 1n \, \E \bigg ( \bigg |  \sum_{k=1}^n (X_k -Y_k) \bigg| \bigg) 
$$
from which the claim follows by convergence of moments in the central 
limit theorem.

When $d \geq 3$, a standard argument (cf.~e.g.~\cite {D69}) goes as follows,
for independent random variables $X_1, \ldots, X_n$ with common
uniform distribution $\mu$ on $[0,1]^d$.
By the Kantorovich-Rubinstein representation \eqref {eq.kr} of $\mathrm {W}_1$,
$$
\mathrm {W}_1 (\mu_n, \mu) \, \geq \, \int_{[0,1]^d} 
\mathrm {dist} \big (x, \{X_1, \ldots, X_n \} \big)\, d\mu(x) .
$$
Let $C_\ell$, $\ell = 1, \ldots, n$, be partition of $[0,1]^d$ into $n$ 
cubes with length $\frac {1}{n^{1/d}}$, so that
$$
\E \big ( \mathrm {W}_1 (\mu_n, \mu) \big)
	\, \geq \, \sum_{\ell = 1}^n \, \E \bigg (\int_{C_\ell} 
	\mathrm {dist} \big (x, \{X_1, \ldots, X_n \} \big)\, d\mu(x) \bigg).
$$
If $D_\ell$ is the collection of cubes surrounding $C_\ell$, then
$\PP ( \forall \, k= 1, \ldots, n ; \,  X_k \notin D_\ell) \geq c $ 
for some $c >0$ only depending on $d$. As a result,
$$
\E \big ( \mathrm {W}_1 (\mu_n, \mu) \big)
	 \, \geq \, \sum_{\ell = 1}^n \, \frac {c}{n^{1/d}} \, \mu (C_\ell) 
	 \, = \, \frac {c}{n^{1/d}} \, .
$$

\section {Quantitative bounds} \label {sec.8}

In this last section, we briefly investigate quantitative bounds, 
in particular with respect to dependence on the dimensional constant $d$
in the main statement (Theorem~\ref {thm.main}) in the form of 
\eqref {eq.quantitative}. We somewhat expand the framework to cover 
at the same time Corollary~\ref {cor.tau}.

\begin {proposition} \label {prop.quantitative}
Let $\mu$ and $\nu$ be two random probability measures on the cube 
$[0,1]^d$ such that their characteristic functions satisfy
$\E( |f_\mu(\pi m) - f_\nu (\pi m) |^2 ) \leq \delta^2$ for all $m \in \Z^d$
for some $0 \leq \delta \leq 2$. Then
\beqs 
\E \big ( \mathrm {W}_1 (\mu, \nu) \big)  \, \leq \,
\begin {cases}
\delta & \text{if \, $d=1$,} \\
5 \delta \, \sqrt { 1 + \log \big ( \frac {4}{\delta^2} \big)} & 
          \text{if \, $d=2$,} \\
10 \sqrt {d} \, \delta^{2/d}   & \text {if \, $ d \geq 3$.} \\
\end {cases}
\eeqs
\end {proposition}

This proposition applied with $\mu=\mu_n$, $\nu = \nu_n$ and 
$\delta = \frac {2}{\sqrt n}$ in the setting of Section~\ref {sec.5} yields 
\eqref {eq.quantitative}. Applied to $\mu = \mu_\tau$,
$\nu = \E_{\pi_n}(\mu_\tau)$ and $\delta = \sqrt {\frac 2n}$
in the setting of Section~\ref {sec.6}, it yields Corollary~\ref {cor.tau}. 

\begin {proof} Using the homogeneity of the distance $\mathrm {W}_1$, 
one may equivalently formulate Proposition~\ref {prop.quantitative} for 
random measures $\mu$ and $\nu$ supported on the cube $[0,\pi]^d$ as
\beqs 
\E \big ( \mathrm {W}_1 (\mu, \nu) \big)  \, \leq \,
\begin {cases}
\pi \delta   & \text{if \, $d=1$,} \\
5\pi\, \delta \, \sqrt { 1 + \log \big ( \frac {4}{\delta^2} \big)} & 
               \text{if \, $d=2$,} \\
10 \pi\, \sqrt {d} \, \delta^{2/d}   & \text {if \, $ d \geq 3$,} \\
\end {cases}
\eeqs
under the assumption that 
$$
\E\big (\big |f_\mu(m) - f_\nu (m) \big|^2 \big) \, \leq \, \delta^2 \quad
	\mbox {for all} \, \,  m \in \Z^d.
$$
That is, the resulting inequalities for measures supported on the 
standard cube $[0,1]^d$ rather than on $[0,\pi]^d$ are obtained with 
numerical factors divided by $\pi$.

Given thus two random measures $\mu$ and $\nu$ on some probability space 
$(\Omega, \mathcal {A}, \PP)$ supported on the cube $[0,\pi]^d$ 
satisfying the latter, averaging the inequality of 
Proposition~\ref {prop.main} yields for every $t >0$,
\beq \label {eq.optimtdelta}
\E \big ({\mathrm {W}}_1(\mu,\nu) \big)
\, \leq \, \delta \, \bigg (\sum_{|m| > 0} \frac {1}{|m|^2} \, 
e^{-2|m|^2t} \bigg)^{1/2} +  2\, \sqrt {2dt} .
\eeq
The task is therefore to suitably optimize in $t>0$.

First note that when $d=1$, the smoothing operation is actually 
not needed and we may simply take $t \to 0$ to get that
$$
\E \big ({\mathrm {W}}_1(\mu,\nu) \big) \, \leq \, 
\frac {\pi}{\sqrt 3} \, \delta \, \leq \, \pi \delta.
$$

Let us then examine more specifically the cases $d=2$ and $d \geq 3$ 
analyzing the sum
$$
\widetilde {S}_d(t) \, = \, S_d\Big(\frac t2\Big) \,= \, 
		\sum_{|m| >0} \frac{1}{|m|^2}\, e^{-|m|^2 t} , \quad t >0.
$$ 
The function $\widetilde {S}_d(t)$ is decreasing in $t>0$, 
vanishing at infinity, and
\beq \label {eq.td}
T_d(t) \, = \, - \widetilde {S}_d'(t) \, = \,  \sum_{|m| > 0} e^{- |m|^2t} 
 		\, = \, \big (1 + T_1(t)\big)^d - 1.
\eeq
In view of the monotonicity of the function $x \to e^{-t x^2}$ for $x>0$, 
we have
$$
\sum_{m=2}^\infty e^{- m^2 t} \, \leq \, \int_1^\infty \! e^{-t x^2}dx
 	\, = \, \frac{1}{\sqrt{2t}} \int_{\sqrt{2t}}^\infty e^{-y^2/2}dy  \,\leq \, 
			\frac{\sqrt \pi}{2 \sqrt t}\ e^{-t}.
$$
Hence, for any $t>0$,
$$
T_1(t) \, = \, \sum_{m \in \Z \setminus {\{0\}}}  e^{-m^2 t} 
	\, = \,  2\,e^{-t} + 2 \sum_{m = 2}^\infty e^{-m^2 t} \, \leq \, 
   		\bigg (2 + \sqrt {\frac{\pi}{t}} \, \bigg)\, e^{-t}.
$$
Putting $a = \sqrt{t}$ and $b = (2\sqrt{t} + \sqrt{\pi} \,)\,e^{-t}$, 
it holds that
\beqs \begin {split}
t^{d/2}\,\big[\big(1 + T_1(t)\big)^d - 1\big] 
 & \, \leq \, (a + b)^d - a^d \\
 & \, = \, \sum_{\ell =0}^{d-1} \binom{d}{\ell} a^\ell \, b^{d-\ell} \\ 
 & \, \leq \,   \big (2 \sqrt{t} + \sqrt{\pi} \, \big)^d \, 
\sum_{\ell =0}^{d-1} \binom{d}{\ell} e^{-(d-\ell)\,t} \\
 &\, \leq \, 2^d\big(2\sqrt{t} + \sqrt{\pi} \, \big )^d \, e^{-t}.
\end {split} \eeqs
Hence, from \eqref {eq.td},
$$
T_d(t) \, \leq \, 2^d\,\bigg(2 + \frac{\sqrt{\pi}}{\sqrt{t}}\,\bigg)^d\, e^{-t}.
$$
It follows that, in the range $t \geq \pi$, $T_d(t) \leq 6^d\, e^{-t}$ 
and thus
$$
\widetilde {S}_d(t) \, = \, \int_t^\infty \! T_d(s)ds \, \leq \,  
6^d\,e^{-t} \, \leq \, 6^d\,e^{-\pi}.
$$
On the other hand, if $t \leq \pi$, then  
$$
T_d(t) \, \leq \, 2^d\bigg (\frac{3\sqrt{\pi}}{\sqrt{t}} \,\bigg)^d e^{-t} 
 	\, \leq \, (6\sqrt{\pi} \, )^d\,t^{-d/2}
$$
so that, in the case $d \geq 3$,
\beqs \begin {split}
\widetilde {S}_d(t) 
 	&  \, = \, \widetilde S_d(\pi) + \int_t^{\pi} T_d(s)\,ds  \\
 	& \, \leq \,  6^d\,e^{-\pi} + 2\,(6\sqrt{\pi}\, )^d\,t^{1-(d/2)}  \\
 	 & \, \leq \,  36\, (e^{-\pi} + 2\pi)\,\Big(\frac{36\,\pi}{t}\Big)^{(d/2) - 1}
\end {split} \eeqs
while for $d=2$,
$$
\widetilde {S}_2(t) \, = \, 
	\widetilde {S}_2(\pi) + \int_t^{\pi} T_2(s)\,ds \, \leq \, 
		36\,e^{-\pi} + 36\pi \log \Big( \frac \pi t \Big).
$$

Let us now return to \eqref {eq.optimtdelta} which states that for any $t>0$,
$$
\E \big ({\mathrm {W}}_1(\mu,\nu) \big)
	\, \leq \, \delta \, \sqrt {\widetilde {S}_d(2t) } +  2\, \sqrt {2dt} \, .
$$
When $d \geq 3$, choose $t = 18\pi \delta^{4/d}$ which is less than or equal 
to $\frac \pi 2$ whenever $\delta ^{2/d} \leq \frac {1}{6}$ in which case
\beqs \begin {split}
\E \big ({\mathrm {W}}_1(\mu,\nu) \big)
 & \, \leq \,   
6\,\Big (\sqrt{ e^{-\pi} + 2\pi} + 2\sqrt { \pi d}\,  \Big)\, \delta^{2/d} \\
 & \, \leq \,  
6\,\Big (\sqrt { e^{-\pi} + 2\pi}\, \sqrt{d/3} + 
2\sqrt { \pi d}\,  \Big)\, \delta^{2/d} \\
  & \, \leq \, 30 \sqrt {d} \, \delta^{2/d}  . \\
\end {split} \eeqs
On the other hand 
$\mathrm {W}_1(\mu,\nu) \leq \pi \sqrt { d}$ for all probability measures 
$\mu$ and $\nu$ supported on $[0,\pi]^d$ so that if 
$\delta ^{2/d} \geq \frac {1}{6}$, the latter inequality is still true.
A similar analysis in the case $d=2$ yields that
$$
\E \big ({\mathrm {W}}_1(\mu,\nu) \big)
		 \, \leq \, 14\, \delta \, 
\sqrt { 1 + \log \Big ( \frac {4}{\delta^2} \Big)} \, .
$$
The proof of the proposition is therefore complete.
\end {proof}

\vskip 8mm

\font\tenrm =cmr10  {\tenrm

\parskip 0mm

\noindent School of Mathematics, University of Minnesota, Minneapolis, MN 55455 USA,
bobkov@math.umn.edu 

\medskip

\noindent Institut de Math\'ematiques de Toulouse, Universit\'e de Toulouse -- Paul-Sabatier, F-31062 Toulouse, France \&  Institut Universitaire de France, ledoux@math.univ-toulouse.fr

}


\vskip 10mm



\begin{thebibliography}{9}


\bibitem {AF84}
E. Acerbi, N. Fusco. Semicontinuity problems in the calculus of variations. 
\textit {Arch. Rational Mech. Anal.}~86, 125--145 (1984).

\bibitem {AKT84}
M. Ajtai, J. Koml\'os, G. Tusn\'ady. On optimal matchings. \textit {Combinatorica}~4, 259--264 (1984).

\bibitem {AST19}
L. Ambrosio, F. Stra, D. Trevisan. A PDE approach to a 2-dimensional matching problem.
\textit {Probab. Theory Related Fields}~173, 433--478 (2019).

\bibitem {BB13}
F. Barthe, C. Bordenave. Combinatorial optimization over two random point sets.
\textit {Séminaire de Probabilités XLV, Lecture Notes in Mathematics}~2078, 483–535. Springer (2013).

\bibitem {BL16}
S. Bobkov, M. Ledoux. One-dimensional empirical measures, order statistics, and Kantorovich
transport distances (2016). To appear in \textit {Memoirs Amer. Math. Soc.}




\bibitem {B05}
R. Bradley. Basic properties of strong mixing conditions. A survey and some open questions. 
Update of, and a supplement to, the 1986 original. \textit {Probab. Surv.}~2, 107--144 (2005).

\bibitem {B07}
R. Bradley. \textit {Introduction to strong mixing conditions,} vol. 1-3.
Kendrick Press (2007).

\bibitem {CS91}
E. Coffman, P. Shor. A simple proof of the $O( \sqrt {n} \log^{3/4} n)$
upright matching bound. \textit {SIAM J. Discrete Math.}~4, 48--57 (1991).

\bibitem {DSS13}
S. Dereich, M. Scheutzow, R. Schottstedt. Constructive quantization: approximation by empirical
measures. \textit {Ann. Inst. Henri Poincaré Probab. Stat.}~49, 1183--1203 (2013).


\bibitem {D69}
R. Dudley. The speed of mean Glivenko-Cantelli convergence.
\textit {Ann. Math. Statist} 40, 40--50 (1969).

\bibitem {D02}
R. Dudley. \textit {Real analysis and probability.} Revised reprint of the 1989 original. 
Cambridge Studies in Advanced Mathematics~74. Cambridge University Press (2002).

\bibitem {FG15}
N. Fournier, A. Guillin. On the rate of convergence in Wasserstein distance
of the empirical measure. \textit {Probab. Theory Related Fields}~162, 707--738 (2015).

\bibitem {HS92}
M. Hahn, Y. Shao. An exposition of Talagrand's mini-course on matching theorems. 
\textit {Probability in Banach spaces}~8, Progr. Probab.~30, 3--38. Birkh\"auser (1992).

\bibitem {HPZ18}
N. Holden, Y. Peres, A. Zhai. Gravitational allocation for uniform points on the sphere.
\textit {Proc. Natl. Acad. Sci. USA}~115, 9666--9671 (2018).


\bibitem {L17}
M. Ledoux. On optimal matching of Gaussian samples. \textit {Zap. Nauchn. Sem. S.-Petersburg. Otdel. Mat. Inst. Steklov. (POMI) 457, Veroyatnost' i Statistika.}~25, 226--264 (2017).

\bibitem {LS89}
T. Leighton, P. Shor. Tight bounds for minimax grid matching with applications
to the average case analysis of algorithms. \textit {Combinatorica}~9, 161--187 (1989).

\bibitem {P18}
R. Peyre. Comparison between $\mathrm{W}_2$ distance and 
$\dot{\mathrm {H}}^{-1}$ norm, and localization
of Wasserstein distance. \textit {ESAIM Control Optim. Calc. Var.}~24, 1489–1501 (2018).

\bibitem {RR98}
S. T. Rachev, L. R\"uschendorf. \textit {Mass transportation problems}, vol. 1 \& 2. Springer (1998).

\bibitem {S15}
F. Santambrogio. \textit {Optimal Transport for Applied Mathematicians.} 
Progress in Nonlinear Differential Equations and Their Applications. Birkh\"auser (2015).

\bibitem {S85}
P. Shor. Random planar matching and bin packing. Ph.D. Thesis, M.I.T. (1985).

\bibitem {S91}
P. Shor. How to pack better than Best Fit: Tight bounds for average-case on-line bin packing.
\textit {Proc. 32nd Annual Symposium on Foundations of Computer Sciences}, 752--759 (1991).

\bibitem {SY91}
P. Shor, J. Yukich. Minimax grid matching and empirical measures.
\textit {Ann. Probab.}~19, 1338--1348 (1991).

\bibitem {S18}
S. Steinerberger. Wasserstein distance, Fourier series and applications (2018).

\bibitem {T92a}
M. Talagrand. Matching random samples in many dimensions. 
\textit {Ann. Appl. Probab.}~2, 846--856 (1992).

\bibitem {T92b}
M. Talagrand. The Ajtai-Komlós-Tusnády matching theorem for general measures.
\textit {Probability in Banach spaces}~8, Progr. Probab.~30, 39--54. Birkh\"auser (1992).

\bibitem {T94a}
M. Talagrand. The transportation cost from the uniform measure to the empirical
measure in dimension $\geq 3$. \textit {Ann. Probab.}~22, 919--959 (1994).

\bibitem {T94b}
M. Talagrand. Matching theorems and empirical discrepancy computations using majorizing measures.
\textit {J. Amer. Math. Soc.}~7, 455--537 (1994).

\bibitem {T14}
M. Talagrand.
\textit {Upper and lower bounds of stochastic processes. Modern methods and classical problems.}
Ergebnisse der Mathematik und ihrer Grenzgebiete~60. Springer (2014).

\bibitem {TY93}
M.~Talagrand, J.~Yukich.
The integrability of the square exponential transportation cost.
\textit {Ann. Appl. Probab.}~3, 1100--1111 (1993).


\bibitem{V03}
C. Villani. \textit {Topics in optimal transportation.} Graduate Studies in 
Mathematics, vol.~58. American Mathematical Society (2003).

\bibitem {Y92}
J. Yukich. Some generalizations of the Euclidean two-sample matching problem. 
\textit {Probability in Banach spaces}~8, Progr. Probab.~30, 55--66. Birkh\"auser (1992).

\bibitem {Y98}
J. Yukich. \textit {Probability theory of classical Euclidean optimization problems.} 
Lecture Notes in Mathematics~1675. Springer (1998).


\end{thebibliography}
\end {document}